\documentclass{article}%
\usepackage{spie}
\usepackage{amsmath}
\usepackage[dvips]{graphicx}%
\usepackage{amsfonts}%
\usepackage{amssymb}
\newtheorem{corollary}[theorem]{Corollary}
\newtheorem{sublemma}[theorem]{Sublemma}
\newtheorem{remark}[theorem]{Remark}
\mathcode`\;="303B
\newenvironment{proofl}
{\par{\it Proof of Lemma \textup{\ref{LemRep.1}}}. \ignorespaces}
{\endproof}

\begin{document}

\title{Compactly supported wavelets and representations of the\\Cuntz relations, II}
\author{Palle E. T. Jorgensen$^{*}$\\Department of Mathematics, The University of Iowa, Iowa City, IA 52242-1419 U.S.A.}
\date{}
\maketitle

\begin{abstract}
We show that compactly
supported wavelets in $L^{2}\left( \mathbb{R}\right) $
of scale $N$ may be effectively
parameterized with a finite
set of spin vectors in $\mathbb{C}^{N}$, and
conversely that every set of spin vectors
corresponds to a wavelet.
The characterization is given
in terms of irreducible
representations of orthogonality
relations defined from
multiresolution wavelet
filters.
\end{abstract}


\section{Introduction}

\label{Int}

\footnotetext{\hskip-1.95\parindent $^{*}$\textbf{E-mail:} \texttt{jorgen@math.uiowa.edu}}Let $L^{2}\left(  \mathbb{R}\right)  $ be the Hilbert space
of all $L^{2}$-functions on $\mathbb{R}$. Let $\psi\in L^{2}\left(
\mathbb{R}\right)  $, and set%
\begin{equation}
\psi_{n,k}\left(  x\right)  :=2^{\frac{n}{2}}\psi\left(  2^{n}x-k\right)
\text{\qquad for }x\in\mathbb{R}\text{, and }n,k\in\mathbb{Z}. \label{eqInt.1}%
\end{equation}
We say that $\psi$ is a wavelet (in the strict sense) if%
\begin{equation}
\left\{  \psi_{n,k};n,k\in\mathbb{Z}\right\}  \label{eqInt.2}%
\end{equation}
constitutes an orthonormal basis in $L^{2}\left(  \mathbb{R}\right)  $; and we
say that $\psi$ is a wavelet in the frame sense (tight frame) if%
\begin{equation}
\left\|  f\right\|  _{L^{2}\left(  \mathbb{R}\right)  }^{2}=\sum
_{n,k\in\mathbb{Z}}\left|  \left\langle \psi_{n,k}\mid
f\right\rangle \right|  ^{2} \label{eqInt.3}%
\end{equation}
holds for all $f\in L^{2}\left(  \mathbb{R}\right)  $, where $\left\langle
\,\cdot\,\mid
\,\cdot\,\right\rangle $ is the
usual $L^{2}\left(  \mathbb{R}\right)  $-inner product, i.e.,%
\begin{equation}
\left\langle \psi_{n,k}\mid
f\right\rangle =\int_{\mathbb{R}}\overline{\psi_{n,k}\left(
x\right)  }f\left(  x\right)  \,dx. \label{eqInt.4}%
\end{equation}
It is known\cite{Dau92,Hor95} that a given wavelet $\psi$ in the
sense of frames is a (strict) wavelet if and only if $\left\|  \psi\right\|
_{L^{2}\left(  \mathbb{R}\right)  }=1$. Either way, the numbers in
(\ref{eqInt.4}) are the wavelet coefficients.

We shall have occasion to consider scaling on $\mathbb{R}$ other than the
dyadic one, say $x\mapsto Nx$ where $N\in\mathbb{N}$, $N>2$. Then the analogue
of (\ref{eqInt.1}) is%
\begin{equation}
\psi_{n,k}\left(  x\right)  :=N^{\frac{n}{2}}\psi\left(  N^{n}x-k\right)
,\qquad x\in\mathbb{R},\;n,k\in\mathbb{Z}. \label{eqInt.5}%
\end{equation}
However, in that case, it is generically not enough to consider only one
$\psi$ in $L^{2}\left(  \mathbb{R}\right)  $: If the wavelet is constructed
from an $N$-subband wavelet filter as in Ref.\ \citenum{BrJo00}, then we will be able
to construct $\psi^{\left(  1\right)  },\psi^{\left(  2\right)  },\dots
,\psi^{\left(  N-1\right)  }$ in $L^{2}\left(  \mathbb{R}\right)  $ such that
the functions in (\ref{eqInt.5}) have the basis property, either in the strict
sense, or in the sense of frames. In that case, the system%
\begin{equation}
\left\{  \psi_{n,k}^{\left(  i\right)  };1\leq i<N,\;n,k\in\mathbb{Z}\right\}
\label{eqInt.6}%
\end{equation}
will constitute an orthonormal basis of $L^{2}\left(  \mathbb{R}\right)  $,
or, alternatively, a tight frame, as in (\ref{eqInt.3}) but with the
$\psi_{n,k}^{\left(  i\right)  }$ functions in place of $\psi_{n,k}$.

In principle, there are many ways (see below) of constructing wavelets
(\ref{eqInt.6}), but we will show in this paper that the method of quadrature
mirror wavelet filters (QMF) has some features that set it apart from the
alternative constructions. Several of the constructions are based on frequency
subbands, and the subbands correspond to a family of closed subspaces in the
Hilbert space $L^{2}\left(  \mathbb{R}\right)  $, but we will show that this
subspace structure is ``optimal'' for the QMF wavelets, in the sense that the
subspaces cannot in a nontrivial way be refined into additional subbands. We
will formulate this result in a mathematically precise fashion, which is based
on a representation of the operators which define the QMF's. In fact, we give
a formula for all these QMF's in the case of compactly supported wavelets.

The present paper is concerned with the wavelet filters which enter into the
construction of the functions $\psi^{(1)},\psi^{(2)},\dots,\psi^{(N-1)}$.
These filters (see (\ref{eqInt.7})--(\ref{eqInt.9}) and (\ref{eq2.15}) below)
are really just a finite set of numbers which relate the $\mathbb{Z}%
$-translates of these functions to the corresponding scalings by $x\mapsto
Nx$. Hence the analysis of the wavelets may be discretized via the filters,
but the question arises whether or not the data which go into the wavelet
filters are \emph{minimal}. It turns out that representation theory is ideally
suited to make the minimality question mathematically precise. (This is a
QMF-multiresolution construction, and it is the minimality and efficiency of
this construction which concern us here. While it is true, see, e.g.,
Refs.\ \citenum{Gab98,FPT99,Bag99,BaMe99,DaLa98},
that there are other and different possible wavelet constructions, it is not
yet clear how our present techniques might adapt to the alternative
constructions, although the approach in Ref.\ \citenum{DaLa98} is also based on
operator-theoretic considerations.)

\section{The scaling function}

\label{Sca}

Since the wavelets come from
the multiresolution functions, it
is of interest to give explicit
constructions for these. We do this
in the present paper, which
introduces two new tools for
explicit constructions of multiresolution
wavelet filters, (i) a certain infinite-dimensional loop
group, and (ii) a certain family
of irreducible representations of orthogonality
relations (the Cuntz relations). Our viewpoint makes
it clear, in particular, that
compactly supported wavelets
may be specified effectively
with a finite set of
spin vectors in $\mathbb{C}^{N}$ ($N$ is the
wavelet scaling number, e.g., $N=2$
in the dyadic case). Hence
these wavelets are given by
a finite set (in arbitrary configuration) of $k$ $Q$-bits
where it turns out that $k$ is
half of the length of the support
of the wavelet in question.

To better explain the minimality issue for multiresolution quadrature mirror
(QMF) wavelet filters, we recall the scaling function $\varphi$ of a
resolution in $L^{2}\left(  \mathbb{R}\right)  $.

Let $g\in\mathbb{N}$, and let $a_{0},a_{1},\dots,a_{2g-1}$ be given complex
numbers such that
\begin{equation}
\sum_{k=0}^{2g-1}a_{k}=2, \label{eqInt.7}%
\end{equation}
and%
\begin{equation}
\sum_{k}a_{k+2l}\bar{a}_{k}=%
\begin{cases}
2 & \text{\quad if }l=0,\\
0 & \text{\quad if }l\neq0.
\end{cases}
\label{eqInt.8}%
\end{equation}
In the formulation of (\ref{eqInt.8}), and elsewhere, we adopt the convention
that terms in a summation are defined to be zero when the index is not in the
specified range: Hence, in (\ref{eqInt.8}), it is understood that $a_{k+2l}=0$
whenever $k$ and $l$ are such that $k+2l$ is not in $\left\{  0,1,\dots
,2g-1\right\}  $. When $\left\{  a_{k};k=0,1,\dots,2g-1\right\}  $ is given
subject to (\ref{eqInt.7})--(\ref{eqInt.8}), then it is known
\cite{BrJo00,BEJ00,Mal99} that there is a $\varphi\in L^{2}\left(
\mathbb{R}\right)  \setminus\left\{  0\right\}  $, unique up to a constant
multiple, such that
\begin{equation}
\varphi\left(  x\right)  =\sum_{k=0}^{2g-1}a_{k}\varphi\left(  2x-k\right)
,\qquad x\in\mathbb{R}, \label{eqInt.9}%
\end{equation}
and $\varphi$ is of compact support; in fact, then
\begin{equation}
\operatorname*{supp}\left(  \varphi\right)  \subset\left[  0,2g-1\right]  .
\label{eqInt.10}%
\end{equation}
(If $H$ denotes the Hilbert transform of $L^{2}\left(  \mathbb{R}\right)  $,
and $\varphi$ solves (\ref{eqInt.9}), then $H\varphi$ does as well; but
$H\varphi$ will not be of compact support if (\ref{eqInt.10}) holds.) In
finding $\varphi$ in (\ref{eqInt.9}), there are methods based on iteration,
on random matrix products, and on the Fourier transform, see
Refs.\ \citenum{Dau92}, \citenum{BrJo00}, \citenum{BEJ00}, and 
\citenum{BrJo99b,Jor00,Coh92b,CoRy95};
and the various methods intertwine in the analysis of $\varphi$,
i.e., in deciding when $\varphi\left(  x\right)  $ is continuous, or not, or
if it is differentiable.

Let $\varphi$ be as in (\ref{eqInt.9}), and let $\mathcal{V}_{0}$ be the
closed subspace in $\mathcal{H}$ ($:=L^{2}\left(  \mathbb{R}\right)  $)
spanned by $\left\{  \varphi\left(  x-k\right)  ;k\in\mathbb{Z}\right\}  $,
i.e., by the integral translates of the scaling function $\varphi$. Let $U$
($:=U_{N}$) be
\begin{equation}
Uf\left(  x\right)  :=N^{-\frac{1}{2}}f\left(  \frac{x}{N}\right)  ,\qquad
f\in L^{2}\left(  \mathbb{R}\right)  , \label{eqInt.11}%
\end{equation}
the unitary scaling operator in $\mathcal{H}=L^{2}\left(  \mathbb{R}\right)
$. Then%
\begin{equation}
U\mathcal{V}_{0}\subset\mathcal{V}_{0}, \label{eqInt.12}%
\end{equation}
it is a proper subspace, and
\begin{equation}
\bigwedge_{n}U^{n}\mathcal{V}_{0}=\left\{  0\right\}  ; \label{eqInt.13}%
\end{equation}
see Ref.\ \citenum{BEJ00} and Ch.~5 of Ref.\ \citenum{Dau92}. 
For $N=2$, the situation is as in
Table \ref{TGA}. Setting
\begin{equation}
\mathcal{V}_{n}:=U^{n}\mathcal{V}_{0}, \label{eqInt.14}%
\end{equation}
and%
\begin{equation}
\mathcal{W}_{n}:=\mathcal{V}_{n-1}\ominus\mathcal{V}_{n}, \label{eqInt.15}%
\end{equation}
we arrive at the resolution%
\begin{equation}
\mathcal{V}_{0}=\sum_{n\geq1} \!\vphantom{\sum}^{\raisebox{3pt}%
[0pt][0pt]{$\scriptstyle\oplus$}}\, \mathcal{W}_{n}, \label{eqInt.16}%
\end{equation}
and the wavelet function $\psi$ is picked in $\mathcal{W}_{0}$. We will set up
an isomorphism between the resolution subspace $\mathcal{V}_{0}$ and $\ell
^{2}\left(  \mathbb{Z}\right)  $, and associate operators in $\ell^{2}\left(
\mathbb{Z}\right)  $ with the wavelet operations in $\mathcal{V}_{0}\subset
L^{2}\left(  \mathbb{R}\right)  $. This is of practical significance given
that the operators in $\ell^{2}\left(  \mathbb{Z}\right)  $ are those which
are defined directly from the wavelet filters, and it is the digital filter
operations which lend themselves to algorithms.

For the general case of scale $N$ ($>2$) the space $\mathcal{V}_{0}\ominus
U_{N}\mathcal{V}_{0}$ splits up as a sum of orthogonal spaces $\mathcal{W}_{1}
^{\left(  i\right)  }$, $i=1,2,\dots,N-1$.

\begin{table}[ptb]
\caption{Discrete vs.\ continuous wavelets, i.e., $\ell\sp2$ vs.\ $L\sp
2\left(  \mathbb{R}\right)  $}%
\label{TGA}
\begin{center}
$\renewcommand{\arraystretch}{0.25}%
\begin{tabular}
[c]{ccccrrrcl}%
$\left\{  0\right\}  $ & $\longleftarrow$ & $\cdots$ & $\longleftarrow$ &
$\mathcal{V}_{2}\raisebox{-12pt}{$\searrow\vphantom{\raisebox{-2pt}{$%
\searrow$}}$\hskip-8pt}$ & $\mathcal{V}_{1}\raisebox{-12pt}{$\searrow
\vphantom{\raisebox{-2pt}{$\searrow$}}$\hskip-8pt}$ & $\mathcal{V}%
_{0}\raisebox{-12pt}{$\searrow\vphantom{\raisebox{-2pt}{$\searrow$}}$%
\hskip-8pt}$ &  & finer scales\\\hline
&  &  &  & \multicolumn{1}{c}{} & \multicolumn{1}{c}{} & \multicolumn{1}{c}{}
&  & \\\cline{1-7}
&  &  &  & \multicolumn{1}{c}{} & \multicolumn{1}{c}{} & \multicolumn{1}{c}{}
& \multicolumn{1}{|c}{} & \\\cline{1-6}
&  &  &  & \multicolumn{1}{c}{} & \multicolumn{1}{c}{} & \multicolumn{1}{|c}{}
& \multicolumn{1}{|c}{} & \\\cline{1-5}
&  &  &  & \multicolumn{1}{c}{} & \multicolumn{1}{|c}{} &
\multicolumn{1}{|c}{} & \multicolumn{1}{|c}{} & \\\cline{1-4}
&  & $\cdots$ &  & \multicolumn{1}{|c}{$\mathcal{W}_{3}$} &
\multicolumn{1}{|c}{$\mathcal{W}_{2}$} & \multicolumn{1}{|c}{$\mathcal{W}_{1}%
$} & \multicolumn{1}{|c}{$\cdots$} & rest of $L^{2}\left(  \mathbb{R}\right)
$\\\cline{1-4}
&  &  &  & \multicolumn{1}{c}{} & \multicolumn{1}{|c}{} &
\multicolumn{1}{|c}{} & \multicolumn{1}{|c}{} & \\\cline{1-5}
&  &  &  & \multicolumn{1}{c}{} & \multicolumn{1}{c}{} & \multicolumn{1}{|c}{}
& \multicolumn{1}{|c}{} & \\\cline{1-6}
&  &  &  & \multicolumn{1}{c}{} & \multicolumn{1}{c}{} & \multicolumn{1}{c}{}
& \multicolumn{1}{|c}{} & \\\cline{1-7}
&  &  &  & \multicolumn{1}{c}{} & \multicolumn{1}{c}{} & \multicolumn{1}{c}{}
&  & \\\hline
\rule[6pt]{0pt}{6pt} &  & $\cdots$ & \multicolumn{1}{r}{$\rlap{$\underset
{\textstyle U}{\longleftarrow}$}\hskip3pt$} & $\rlap{$\underset{\textstyle
U}{\longleftarrow}$}\hskip3pt$ & $\rlap{$\underset{\textstyle U}%
{\longleftarrow}$}\hskip3pt$ & \multicolumn{1}{c}{} &  & \\
&  & $\llap{$W$}\makebox[6pt]{\raisebox{4pt}{\makebox[0pt]{\hss$\uparrow$\hss
}}\raisebox{-4pt}{\makebox[0pt]{\hss$|$\hss}}}$ &  &
\multicolumn{1}{c}{$\makebox[6pt]{\raisebox{4pt}{\makebox[0pt]{\hss$\uparrow
$\hss}}\raisebox{-4pt}{\makebox[0pt]{\hss$|$\hss}}}$} &
\multicolumn{1}{c}{$\makebox[6pt]{\raisebox{4pt}{\makebox[0pt]{\hss$\uparrow
$\hss}}\raisebox{-4pt}{\makebox[0pt]{\hss$|$\hss}}}$} &
\multicolumn{1}{c}{$\makebox[6pt]{\raisebox{4pt}{\makebox[0pt]{\hss$\uparrow
$\hss}}\raisebox{-4pt}{\makebox[0pt]{\hss$|$\hss}}}\rlap{$W$}$} &  & \\
\rule[-6pt]{0pt}{6pt}$\left\{  0\right\}  $ & $\longleftarrow$ & $\cdots$ &
\multicolumn{1}{r}{$\rlap{$\overset{\textstyle S_{0}}{\longleftarrow}$}%
\hskip3pt$} & $\rlap{$\overset{\textstyle S_{0}}{\longleftarrow}$}\hskip3pt$ &
$\rlap{$\overset{\textstyle S_{0}}{\longleftarrow}$}\hskip3pt$ &
\multicolumn{1}{c}{} &  & \\\cline{1-7}
&  &  &  & \multicolumn{1}{c}{} & \multicolumn{1}{c}{} & \multicolumn{1}{c}{}
& \multicolumn{1}{|c}{} & \\\cline{1-6}
&  &  &  & \multicolumn{1}{c}{} & \multicolumn{1}{c}{} & \multicolumn{1}{|c}{}
& \multicolumn{1}{|c}{} & \\\cline{1-5}
&  &  &  & \multicolumn{1}{c}{} & \multicolumn{1}{|c}{} &
\multicolumn{1}{|c}{} & \multicolumn{1}{|c}{} & \\\cline{1-4}
&  & $\cdots$ &  & \multicolumn{1}{|c}{$S_{0}^{2}\mathcal{L}$} &
\multicolumn{1}{|c}{$S_{0}\mathcal{L}$} & \multicolumn{1}{|c}{$\mathcal{L}%
=S_{1}\ell^{2}$} & \multicolumn{1}{|c}{} & \\\cline{1-4}
&  &  &  & \multicolumn{1}{c}{} & \multicolumn{1}{|c}{} &
\multicolumn{1}{|c}{} & \multicolumn{1}{|c}{} & \\\cline{1-5}
&  &  &  & \multicolumn{1}{c}{} & \multicolumn{1}{c}{} & \multicolumn{1}{|c}{}
& \multicolumn{1}{|c}{} & \\\cline{1-6}
&  &  &  & \multicolumn{1}{c}{} & \multicolumn{1}{c}{} & \multicolumn{1}{c}{}
& \multicolumn{1}{|c}{} & \\\cline{1-7}
&  &  &  & \makebox[36pt]{\hfill}\llap{$S_{0}^{2}\ell^{2}$}\raisebox
{12pt}{$\nearrow\vphantom{\raisebox{2pt}{$\nearrow$}}$\hskip-8pt} &
\makebox[36pt]{\hfill}\llap{$S_{0}\ell^{2}$}\raisebox{12pt}{$\nearrow
\vphantom{\raisebox{2pt}{$\nearrow$}}$\hskip-8pt} & $\ell^{2}\raisebox
{12pt}{$\nearrow\vphantom{\raisebox{2pt}{$\nearrow$}}$\hskip-8pt}$ &  &
\end{tabular}
\ \ $
\end{center}
\end{table}

\section{Representations of $\mathcal{O}_{N}$ and Table \ref{TGA}%
\protect\linebreak \ (Discrete vs.\ continuous wavelets)}

\label{Rep}

The practical significance of the operator system in Table \ref{TGA} (scale
$N=2$) is that the operators which generate wavelets in $L^{2}\left(
\mathbb{R}\right)  $ become modeled by an associated system of operators in
the \emph{sequence space} $\ell^{2}$ ($:=\ell^{2}\left(  \mathbb{Z}\right)
\cong L^{2}\left(  \mathbb{T}\right)  $). (We will do the discussion here in
Section \ref{Rep} just for $N=2$, but this is merely for simplicity. It easily
generalizes to arbitrary $N$.) Then the algorithms are implemented in
$\ell^{2}$ by basic discrete operations, and only in the end are the results
then ``translated'' back to the space $L^{2}(\mathbb{R})$. The space
$L^{2}(\mathbb{R})$ is not amenable (in its own right) to \emph{discrete}
computations. This is made precise by virtue of the frame operator
$W\colon\ell^{2}$ ($\cong L^{2}\left(  \mathbb{T}\right)  $) $\rightarrow
\mathcal{V}_{0}$ ($\subset L^{2}\left(  \mathbb{R}\right)  $) which may be
defined as%
\begin{equation}
W\colon\ell^{2}\ni\left(  \xi_{k}\right)  \longmapsto\sum_{k\in\mathbb{Z}}%
\xi_{k}\varphi\left(  x-k\right)  \in L^{2}\left(  \mathbb{R}\right)  .
\label{eqRep.1}%
\end{equation}
If the scaling function $\varphi$ has then been constructed to have orthogonal
translates, then $W$ will be an isometry of $\ell^{2}$ onto $\mathcal{V}_{0}$
($\subset L^{2}\left(  \mathbb{R}\right)  $). Even if the functions $\left\{
\varphi\left(  x-k\right)  \right\}  _{k\in\mathbb{Z}}$ formed from $\varphi$
by $\mathbb{Z}$-translates only constitute a frame in $\mathcal{V}_{0}$, then
we will have the following estimates:%
\begin{equation}
c_{1}\sum_{k\in\mathbb{Z}}\left|  \xi_{k}\right|  ^{2}\leq\int_{-\infty
}^{\infty}\left|  \left(  W\xi\right)  \left(  x\right)  \right|  ^{2}\,dx\leq
c_{2}\sum_{k\in\mathbb{Z}}\left|  \xi_{k}\right|  ^{2}, \label{eqRep.2}%
\end{equation}
where $c_{1}$ and $c_{2}$ are positive constants depending only on the scaling
function $\varphi$, or equivalently,%
\begin{equation}
c_{1}^{1/2}\cdot\left\|  \xi\right\|  _{\ell^{2}}\leq\left\|  W\xi\right\|
_{L^{2}\left(  \mathbb{R}\right)  }\leq c_{2}^{1/2}\cdot\left\|  \xi\right\|
_{\ell^{2}}\,. \label{eqRep.3}%
\end{equation}

\begin{lemma}
\label{LemRep.1}If the coefficients $\left\{  a_{k};k=0,1,\dots,2g-1\right\}
$ from \textup{(\ref{eqInt.9})} satisfy the conditions in
\textup{(\ref{eqInt.8}),} then the corresponding operator $S_{0}\colon\ell
^{2}\rightarrow\ell^{2}$, given by%
\begin{equation}
\left(  S_{0}\xi\right)  _{k}=\frac{1}{\sqrt{2}}\sum_{l\in\mathbb{Z}}%
a_{k-2l}\xi_{l}=\frac{1}{\sqrt{2}}\sum_{\substack{p\in\mathbb{Z}%
\colon\\p\equiv k\operatorname{mod}2}}a_{p}\xi_{\frac{k-p}{2}},\qquad
k\in\mathbb{Z}, \label{eqRep.4}%
\end{equation}
is \emph{isometric} and satisfies the following \emph{intertwining identity:}%
\begin{equation}
WS_{0}=UW, \label{eqRep.5}%
\end{equation}
where $U$ is the dyadic scaling operator in $L^{2}\left(  \mathbb{R}\right)  $
which we introduced in \textup{(\ref{eqInt.5}).} \textup{(}Here we restrict
attention to $N=2$, but just for notational simplicity!\/\textup{)} Setting%
\begin{equation}
b_{k}:=\left(  -1\right)  ^{k}\bar{a}_{2g-1-k}\,, \label{eqRep.6}%
\end{equation}
and defining a second isometric operator $S_{1}\colon\ell^{2}\rightarrow
\ell^{2}$ by formula \textup{(\ref{eqRep.4})} with the only modification that
$\left(  b_{k}\right)  $ is used in place of $\left(  a_{k}\right)  $, we get%
\begin{equation}
S_{j}^{\ast}S_{k}^{{}}=\delta_{j,k}\hbox{\upshape{\small1}\kern-3.3pt1}%
_{\ell^{2}} \label{eqRep.7}%
\end{equation}
and%
\begin{equation}
\sum_{j}S_{j}^{{}}S_{j}^{\ast}=\hbox{\upshape{\small1}\kern-3.3pt1}_{\ell^{2}%
}\,, \label{eqRep.8}%
\end{equation}
which are the Cuntz identities from operator theory\cite{Cun77}, and the
operators $S_{0}$ and $S_{1}$ satisfy the identities indicated in Table
\textup{\ref{TGA}.}
\end{lemma}

\begin{remark}
\label{RemRepNew.pound}For understanding the second line in Table
\textup{\ref{TGA},} note that $S_{0}$ is a \emph{shift} as an isometry, in the
sense of Ref.\ \textup{\citenum{SzFo70},} 
and $\mathcal{L}:=S_{1}\ell^{2}$ is a wandering
subspace for $S_{0}$, in the sense that the spaces $\mathcal{L}$, $S_{0}^{{}%
}\mathcal{L}$, $S_{0}^{2}\mathcal{L},\;\dots$ are mutually orthogonal in
$\ell^{2}$. To see this, note that \textup{(\ref{eqRep.8})} implies that%
\begin{equation}
\left(  \mathcal{L}:=\right)  \qquad S_{1}\ell^{2}=\ell^{2}\ominus S_{0}%
\ell^{2}=\ker\left(  S_{0}^{\ast}\right)  . \label{eqRemRepNew.pound.1}%
\end{equation}

As a result, we get the following:
\end{remark}

\begin{corollary}
\label{corollary2new}The projections onto the orthogonal subspaces in the
second line of Table $1$ which correspond to the $\mathcal{W}_{1}%
,\mathcal{W}_{2},\dots$ subspaces of the first line \textup{(}see
\textup{(\ref{eqInt.15}))} are as follows:
\begin{align*}
\operatorname*{proj}\mathcal{L}  &  =S_{1}S_{1}^{\ast}=I-S_{0}S_{0}^{\ast},\\
\vdots &  \qquad\vdots\\
\operatorname*{proj}S_{0}^{n-1}\mathcal{L}  &  =S_{0}^{n-1}S_{0}^{\ast
n-1}-S_{0}^{n}S_{0}^{\ast n}.
\end{align*}
\end{corollary}

\begin{proof}
Immediate from Lemma \ref{LemRep.1} and Remark \ref{RemRepNew.pound}.
\end{proof}

\begin{remark}
\label{RemRep.2}Any system of operators $\left\{  S_{j}\right\}  $ satisfying
\textup{(\ref{eqRep.7})--(\ref{eqRep.8})} is said to be a representation of
the $C^{\ast}$-algebra $\mathcal{O}_{2}$, and there is a similar notion for
$\mathcal{O}_{N}$ when $N>2$, with $\mathcal{O}_{N}$ having generators
$S_{0},S_{1},\dots,S_{N-1}$, but otherwise also satisfying the operator
identities \textup{(\ref{eqRep.7})--(\ref{eqRep.8}).}
\end{remark}

\begin{definition}
\label{DefRep.3}A representation of $\mathcal{O}_{N}$ on the Hilbert space
$\ell^{2}$ is said to be \emph{irreducible} if there are no closed subspaces
$\left\{  0\right\}  \subsetneqq\mathcal{H}_{0}\subsetneqq\ell^{2}$ which
reduce the representation, i.e., which yield a representation of
\textup{(\ref{eqRep.7})--(\ref{eqRep.8})} on each of the two subspaces in the
decomposition%
\begin{equation}
\ell^{2}=\mathcal{H}_{0}\oplus\left(  \ell^{2}\ominus\mathcal{H}_{0}\right)  ,
\label{eqRep.9}%
\end{equation}
where%
\begin{equation}
\ell^{2}\ominus\mathcal{H}_{0}=\left(  \mathcal{H}_{0}\right)  ^{\perp
}=\left\{  \xi\in\ell^{2};\left\langle \xi\mid
\eta\right\rangle =0,\;\forall
\,\eta\in\mathcal{H}_{0}\right\}  . \label{eqRep.10}%
\end{equation}
\end{definition}

\begin{proofl}
Most of the details of the proof are
contained in Refs.\ \citenum{BrJo97b} and \citenum{BrJo00}, 
so we only sketch points not
already covered there. The essential step (for the present applications) is
the formula (\ref{eqRep.5}), which shows that $W$ intertwines the isometry
$S_{0}\,$with the restriction of the unitary operator $U\colon f\mapsto
\frac{1}{\sqrt{2}}f\left(  x/2\right)  $ to the resolution subspace
$\mathcal{V}_{0}\subset L^{2}\left(  \mathbb{R}\right)  $. We have:%
\begin{align*}
\left(  UW\xi\right)  \left(  x\right)   &  =\frac{1}{\sqrt{2}}\left(
W\xi\right)  \left(  \frac{x}{2}\right)  &  & \\
&  =\frac{1}{\sqrt{2}}\sum_{k\in\mathbb{Z}}\xi_{k}\varphi\left(  \frac{x}%
{2}-k\right)  &  &  \text{\qquad(by (\ref{eqRep.1}))}\\
&  =\frac{1}{\sqrt{2}}\sum_{k\in\mathbb{Z}}\sum_{l\in\mathbb{Z}}\xi_{k}%
a_{l}\varphi\left(  x-2k-l\right)  &  &  \text{\qquad(by (\ref{eqInt.9}))}\\
&  =\frac{1}{\sqrt{2}}\sum_{p\in\mathbb{Z}}\left(  \sum_{k\in\mathbb{Z}}%
\xi_{k}a_{p-2k}\right)  \varphi\left(  x-p\right)  &  & \\
&  =\sum_{p\in\mathbb{Z}}\left(  S_{0}\xi\right)  _{p}\varphi\left(
x-p\right)  &  &  \text{\qquad(by (\ref{eqRep.4}))}\\
&  =\left(  WS_{0}\xi\right)  \left(  x\right)  &  &  \text{\qquad(by
(\ref{eqRep.1}))}%
\end{align*}
for all $\xi\in\ell^{2}$, and all $x\in\mathbb{R}$. This proves (\ref{eqRep.5}).

The rest of the proof will be given in a form slightly more general than
needed. For later use, we record the following table of operators on the
respective Hilbert spaces $L^{2}\left(  \mathbb{T}\right)  \cong\ell^{2}$ and
$L^{2}\left(  \mathbb{R}\right)  $, and the corresponding transformation rules
with respect to the operator $W$. Let $N$ be the scale number, and let
$\left(  a_{k}\right)  _{k=0}^{Ng-1}$ be given satisfying%
\begin{equation}
\sum_{k\in\mathbb{Z}}a_{k+Nl}\bar{a}_{k}=\delta_{0,l}N \label{eqRep.18}%
\end{equation}
and set $m_{0}\left(  z\right)  :=\frac{1}{\sqrt{N}}\sum_{k=0}^{Ng-1}%
a_{k}z^{k}$, $z\in\mathbb{T}$. Then we have the following transformation
rules.%
\begin{equation}%
\begin{array}
[c]{llll}
& \quad\;\text{\textbf{SCALING}} & \quad\;\text{\textbf{TRANSLATION}} & \\
\displaystyle\vphantom{\frac{1}{N}}L^{2}\left(  \mathbb{R}\right)  \colon &
\displaystyle\quad F\longmapsto\smash{\frac{1}{\sqrt{N}}}F\left(  \frac{x}%
{N}\right)  & \displaystyle\quad F\left(  x\right)  \longmapsto F\left(
x-1\right)  & \text{\quad real wavelets}\\
\displaystyle\uparrow\rlap{$\scriptstyle W$} &  &  & \\
\displaystyle\vphantom{\frac{1}{N}}\ell^{2}\colon & \displaystyle\quad
\xi\longmapsto\smash{\sum_l}a_{k-Nl}\xi_{l} & \displaystyle\quad\left(
\xi_{k}\right)  \longmapsto\left(  \xi_{k-1}\right)  & \text{\quad discrete
model}\\
\displaystyle\uparrow\rlap{\small Fourier transform} &  &  & \\
\displaystyle\vphantom{\frac{1}{N}}L^{2}\left(  \mathbb{T}\right)  \colon &
\displaystyle\quad f\longmapsto m_{0}\left(  z\right)  f\left(  z^{N}\right)
& \displaystyle\quad f\left(  z\right)  \longmapsto zf\left(  z\right)  &
{\normalsize \quad\setlength{\arraycolsep}{0pt}%
\begin{array}
[t]{l}%
\text{periodic model,}\\
\quad\mathbb{T}=\mathbb{R}/2\pi\mathbb{Z}%
\end{array}
}%
\end{array}
\label{eqRep.19}%
\end{equation}
Now the proof may be completed by use of the following sublemma, which we
state just for $N=2$.

\begin{sublemma}
\label{SublemRep.4}Let $\mathfrak{G}$ be a subspace of $L^{2}\left(
\mathbb{T}\right)  $ which is invariant under multiplication by $z^{2}$. Then
there is an $m_{1}\in L^{\infty}\left(  \mathbb{T}\right)  $ such that
$\mathfrak{G}=\left\{  m_{1}\left(  z\right)  f\left(  z^{2}\right)  ;f\in
L^{2}\left(  \mathbb{T}\right)  \right\}  $ and $\frac{1}{2}\sum_{w^{2}%
=z}\left|  m_{1}\left(  w\right)  \right|  ^{2}\equiv1$, $\mathrm{a.a.}%
\,z\in\mathbb{T}$.
\end{sublemma}

\begin{proof}
The proof follows from the Beurling-Lax-Halmos theorem\cite{SzFo70}.
\end{proof}

\textit{Completion of proof of Lemma }\ref{LemRep.1}.\textup{ }Let
$m_{0}\left(  z\right)  :=\frac{1}{\sqrt{2}}\sum_{k}a_{k}z^{k}$. Then we saw
that $S_{0}f\left(  z\right)  :=m_{0}\left(  z\right)  f\left(  z^{2}\right)
$ is isometric in $L^{2}\left(  \mathbb{T}\right)  $, and the complementary
space
\[
\mathfrak{G}:=L^{2}\left(  \mathbb{T}\right)  \ominus S_{0}L^{2}\left(
\mathbb{T}\right)  =\ker\left(  S_{0}^{\ast}\right)  =\left\{  f\in
L^{2}\left(  \mathbb{T}\right)  ;\sum_{w^{2}=z}\overline{m_{0}\left(
w\right)  }f\left(  w\right)  =0\right\}
\]
then satisfies the condition in Sublemma \ref{SublemRep.4}. Let $m_{1}$ be the
function determined from the sublemma. Then, after multiplication by a
suitable $z^{2p}$, we will get $m_{1}\left(  z\right)  =\frac{1}{\sqrt{2}}%
\sum_{k}b_{k}z^{k}$, with the coefficients$\mathcal{\ }b_{k}$ as in
(\ref{eqRep.6}). Setting $S_{1}f\left(  z\right)  :=m_{1}\left(  z\right)
f\left(  z^{2}\right)  $, it follows then that (\ref{eqRep.7})--(\ref{eqRep.8}%
) are satisfied.
\end{proofl}

\begin{remark}
\label{RemRep.4}The significance of irreducibility \textup{(}when
satisfied\/\textup{)} is that the \emph{wavelet subbands} which are indicated
in Table \textup{\ref{TGA}} are then the \emph{only subbands} of the
corresponding multiresolution. We will show that in fact irreducibility holds
\emph{generically,} but it does not hold, for example, for the Haar wavelets.
In the simplest case, the Haar wavelet has $g=2=N$ and the numbers
\begin{equation}%
\begin{pmatrix}
a_{0} & a_{1}\\
b_{0} & b_{1}%
\end{pmatrix}
=%
\begin{pmatrix}
1 & 1\\
1 & -1
\end{pmatrix}
. \label{eqRep.11}%
\end{equation}
Hence, for this representation of $\mathcal{O}_{2}$ on $\ell^{2}$, we may take
$\mathcal{H}_{0}=\ell^{2}\left(  0,1,2,\dots\right)  $, and therefore
$\mathcal{H}_{0}^{\perp}=\ell^{2}\left(  \dots,-3,-2,-1\right)  $. Returning
to the multiresolution diagram in Table \textup{\ref{TGA},} this means that we
get additional subspaces of $L^{2}\left(  \mathbb{R}\right)  $, on top of the
standard ones which are listed in Table \textup{\ref{TGA}.} Specifically, in
addition to%
\begin{equation}
\mathcal{V}_{n}=U^{n}\mathcal{V}_{0}=WS_{0}^{n}\ell^{2} \label{eqRep.12}%
\end{equation}
and%
\begin{equation}
\mathcal{W}_{n}=\mathcal{V}_{n-1}\ominus\mathcal{V}_{n}=WS_{0}^{n-1}S_{1}^{{}%
}\ell^{2}, \label{eqRep.13}%
\end{equation}
we get a new system with ``twice as many'', as follows: $\mathcal{V}%
_{n}^{\left(  \pm\right)  }$ and $\mathcal{W}_{n}^{\left(  \pm\right)  }$,
where%
\begin{align}
\mathcal{V}_{n}^{\left(  +\right)  }  &  =WS_{0}^{n}\left(  \mathcal{H}%
_{0}\right)  ,\label{eqRep.14a}\\
\mathcal{W}_{n}^{\left(  +\right)  }  &  =WS_{0}^{n-1}S_{1}^{{}}\left(
\mathcal{H}_{0}\right)  ; \label{eqRep14b}%
\end{align}
and%
\begin{align}
\mathcal{V}_{n}^{\left(  -\right)  }  &  =WS_{0}^{n}\left(  \mathcal{H}%
_{0}^{\perp}\right)  ,\label{eqRep.15a}\\
\mathcal{W}_{n}^{\left(  -\right)  }  &  =WS_{0}^{n-1}S_{1}^{{}}\left(
\mathcal{H}_{0}^{\perp}\right)  . \label{eqRep.15b}%
\end{align}
For the case of the Haar wavelet, see \textup{(\ref{eqRep.11}),}%
\[
\mathcal{V}_{0}^{\left(  +\right)  }\subset L^{2}\left(  0,\infty\right)
,\qquad\mathcal{V}_{0}^{\left(  -\right)  }\subset L^{2}\left(  -\infty
,0\right)  ,
\]
or rather, $\mathcal{V}_{0}$ consists of finite linear combinations of
$\mathbb{Z}$-translates of
\begin{equation}
\varphi\left(  x\right)  =%
\begin{cases}
1 & \text{ if }0\leq x<1,\\
0 & \text{ if }x\in\mathbb{R}\setminus\left[  0,1\right)  ,
\end{cases}
\label{eqRep.16}%
\end{equation}
alias the step functions of step-size one, i.e., functions in $L^{2}\left(
\mathbb{R}\right)  $ which are constant between $n$ and $n+1$ for all
$n\in\mathbb{Z}$; and%
\begin{equation}
\mathcal{V}_{0}^{\left(  +\right)  }=\mathcal{V}_{0}^{{}}\cap L^{2}\left(
0,\infty\right)  ,\qquad\mathcal{V}_{0}^{\left(  -\right)  }=\mathcal{V}%
_{0}^{{}}\cap L^{2}\left(  -\infty,0\right)  . \label{eqRep.17}%
\end{equation}
Hence we get two separate wavelets, but with translations built on $\left\{
0,1,2,\dots\right\}  $ and $\left\{  \dots,-3,-2,-1\right\}  $. 
In view of the
graphics
of the cascade approximation
to the scaling function \textup{(}see
Refs.\ \textup{\citenum{BrJo99b}} and \textup{\citenum{Jor00}),}
it is perhaps surprising that other wavelets
\textup{(}different from the Haar wavelets\/\textup{)} do not have the
corresponding additional ``positive vs.\ negative'' splitting into subbands
within the Hilbert space $L^{2}\left(  \mathbb{R}\right)  $.
\end{remark}

\begin{remark}
\label{RemRep.5}There are other dyadic Haar wavelets, in addition to
\textup{(\ref{eqRep.16}).} For example, let%
\begin{equation}
\varphi_{k}\left(  x\right)  =%
\begin{cases}
\frac{1}{\sqrt{2k+1}} & \text{ if }0\leq x<2k+1,\\
0 & \text{ if }x\in\mathbb{R}\setminus\left[  0,2k+1\right)  .
\end{cases}
\label{eqRep.17bis}%
\end{equation}
Then it follows that there is a splitting of $\mathcal{V}_{0}$ into orthogonal
subspaces which is analogous to \textup{(\ref{eqRep.17}),} but it has many
more subbands than the two, ``positive vs.\ negative'', which are given in
\textup{(\ref{eqRep.17}),} and which are special to the standard Haar wavelet
\textup{(\ref{eqRep.16}).} For details on these other Haar wavelets, and their
decompositions, we refer the reader to 
Proposition \textup{8.2} of Ref.\ \textup{\citenum{BrJo99a}.} The
corresponding $m$-functions of \textup{(\ref{eqRep.17bis})} are%
\begin{equation}
m_{0}\left(  z\right)  =\frac{1}{\sqrt{2}}\left(  1+z^{2k+1}\right)  ,\qquad
m_{1}\left(  z\right)  =\frac{1}{\sqrt{2}}\left(  1-z^{2k+1}\right)  ,\qquad
z\in\mathbb{T}. \label{eqRep.17ter}%
\end{equation}
Hence, after adjusting the $\mathcal{O}_{2}$-representation $T$ with a
\textup{(}rotation\/\textup{)} $V\in\mathrm{U}_{2}\left(  \mathbb{C}\right)
$, we have%
\begin{equation}
T_{0}f\left(  z\right)  =f\left(  z^{2}\right)  ,\qquad T_{1}f\left(
z\right)  =z^{2k+1}f\left(  z^{2}\right)  ,\qquad f\in L^{2}\left(
\mathbb{T}\right)  \cong\ell^{2}, \label{eqRep.17tetra}%
\end{equation}
and, of course, the two new operators $T_{0},T_{1}$ will satisfy the
$\mathcal{O}_{2}$-identities 
\textup{(\ref{eqRep.7})--(\ref{eqRep.8}),} and the
corresponding representation will have the same reducing subspaces as the one
defined directly from $m_{0}$ and $m_{1}$. The explicit decomposition of the
multiresolution subspaces corresponding to \textup{(\ref{eqRep.17})} may be
derived, via $W$ in Table \textup{\ref{TGA},} from the corresponding
decomposition into sums of irreducibles for the $\mathcal{O}_{2}%
$-representation on $\ell^{2}$ which corresponds to \textup{(\ref{eqRep.17}).}
This means that the corresponding \textup{(\ref{eqRep.9})} which is associated
with \textup{(\ref{eqRep.17bis})} and \textup{(\ref{eqRep.17tetra})} has
\emph{more than two} terms in its subspace configuration.
\end{remark}

\section{Wavelet filters and subbands}

\label{Wav}

The power and the usefulness of the multiresolution subband filters for the
analysis of wavelets and their algorithms was first demonstrated forcefully in
Refs.\ \citenum{CoWi93} and \citenum{Wic93}; 
see especially p.~140 of Ref.\ \citenum{CoWi93} and
p.~157 of Ref.\ \citenum{Wic93}, 
where the $\mathcal{O}_{N}$-relations (\ref{eqRep.7}%
)--(\ref{eqRep.8}) are identified, and analyzed in the case $N=2$. Around the
same time, A. Cohen\cite{Coh92b} identified and utilized the interplay
between $\ell^{2}$ and $L^{2}\left(  \mathbb{R}\right)  $ which, as noted in
Section \ref{Rep} above, is implied by the $\mathcal{O}_{N}$-relations and
their representations. But neither of those prior references takes up the
construction of $\mathcal{O}_{N}$-representations in a systematic fashion.

Of course the quadrature mirror filters (QMF's) have a long history in
electrical engineering (speech coding problems), going back to long before
they were used in wavelets, but the form in which we shall use them here is
well articulated, for example, in Ref.\ \citenum{CEG77}. 
Some more of the history of
and literature on wavelet filters is covered well in 
Refs.\ \citenum{Mey93} and
\citenum{Ben00}.

The operators corresponding to wavelet filters may be realized on either one
of the two Hilbert spaces $\ell^{2}(\mathbb{Z})$ or $L^{2}(\mathbb{T})$,
$\mathbb{T}=\mathbb{R}/2\pi\mathbb{Z}$, and $L^{2}(\mathbb{T})$ defined from
the normalized Haar measure $\mu$ on $\mathbb{T}$. But, of course, $\ell
^{2}(\mathbb{Z})$ $\cong L^{2}(\mathbb{T})$ via the Fourier series%

\begin{equation}
f(z)=\sum_{k\in\mathbb{Z}}\xi_{k}z^{k}, \label{eq2.1}%
\end{equation}
and its inverse%
\begin{equation}
\xi_{k}=\hat{f}(k)=\int_{\mathbb{T}}z^{-k}f(z)\,d\mu(z). \label{eq2.2}%
\end{equation}
For a given sequence $a_{0},a_{1},\dots,a_{Ng-1}$, consider the operator
$S_{0}$ in $\ell^{2}(\mathbb{Z})$ given by%

\begin{equation}
\xi\longmapsto S_{0}\xi, \label{eq2.3}%
\end{equation}
and%
\begin{equation}
(S_{0}\xi)_{k}=\frac{1}{\sqrt{N}}\sum_{l}a_{k-lN}\xi_{l}. \label{eq2.4}%
\end{equation}
Setting
\begin{equation}
m_{0}(z)=\frac{1}{\sqrt{N}}\sum_{k=0}^{Ng-1}a_{k}z^{k},\text{ } \label{eq2.5}%
\end{equation}
and%
\begin{equation}
(\hat{S}_{0}f)(z)=m_{0}(z)f(z^{N}),\qquad f\in L^{2}(\mathbb{T}),\qquad
z\in\mathbb{T}, \label{eq2.6}%
\end{equation}
we note that $S_{0}$ and $\hat{S}_{0}$ are really two versions of the same
operator, i.e.,
\begin{equation}
(\hat{S}_{0}f)\widehat{\phantom{S}}=S_{0}(\hat{f}) \label{eq2.7}%
\end{equation}
when $\hat{f}=(\xi_{k})$ from (\ref{eq2.2}). (The first one is the discrete
model, and the second, the periodic model, referring to the diagram
(\ref{eqRep.19}).) Hence, we shall simply use the same notation $S_{0}$ in
referring to this operator in either one of its incarnations. It is the
(\ref{eq2.4}\textbf{) }version which is used in algorithms, of course.

Let $\varphi\in L^{2}(\mathbb{R})$ be the compactly supported scaling function
solving
\begin{equation}
\varphi(x)=\sum_{k=0}^{Ng-1}a_{k}\varphi(Nx-k). \label{eq2.8}%
\end{equation}
Then define the operator
\begin{equation}
W\colon\ell^{2}(\mathbb{Z})\longrightarrow L^{2}(\mathbb{R}) \label{eq2.9}%
\end{equation}
by
\begin{equation}
W\xi=\sum_{k}\xi_{k}\varphi(x-k). \label{eq2.10}%
\end{equation}
The conditions on the wavelet filter $\{a_{k}\}$ from Section \ref{Int} and in
(\ref{eqRep.18}\textbf{), }may now be restated in terms of $m_{0}(z)$ in
(\ref{eq2.5}) as follows:
\begin{equation}
\sum_{k=0}^{N-1}\left|  m_{0}(ze^{i\frac{k2\pi}{N}})\right|  ^{2}=N,\text{ }
\label{eq2.11}%
\end{equation}
and%
\begin{equation}
m_{0}(1)=\sqrt{N}. \label{eq2.12}%
\end{equation}
It then follows from Lemma \ref{LemRep.1} that $W$ in (\ref{eq2.10}) maps
$\ell^{2}(\mathbb{Z})$ onto the resolution subspace $\mathcal{V}_{0}$($\subset
L^{2}(\mathbb{R})$), and that
\begin{equation}
U_{N}W=WS_{0} \label{eq2.13}%
\end{equation}
where
\begin{equation}
U_{N}f(x)=N^{-\frac{1}{2}}f\left(  \frac{x}{N}\right)  ,\qquad f\in
L^{2}(\mathbb{R}),\qquad x\in\mathbb{R}. \label{eq.2.14}%
\end{equation}
We showed in Ref.\ \citenum{BrJo00} 
that there are functions $m_{1},\dots,m_{N-1}$ such
that the $N$-by-$N$ complex matrix
\begin{equation}
\frac{1}{\sqrt{N}}\left(  m_{j}(e^{i\frac{k2\pi}{N}}z)\right)  _{j,k=0}^{N-1}
\label{eq2.15}%
\end{equation}
is unitary for all $z\in\mathbb{T}$. If we define
\begin{equation}
S_{j}f(z)=m_{j}(z)f(z^{N}),\qquad f\in L^{2}(\mathbb{T}),\qquad z\in
\mathbb{T}\text{, } \label{eq2.16}%
\end{equation}
then
\begin{equation}
S_{j}^{\ast}S_{k}=\delta_{j,k}I_{L^{2}(\mathbb{T})}, \label{eq2.17}%
\end{equation}
and%
\begin{equation}
\sum_{j=0}^{N-1}S_{j}S_{j}^{\ast}=I_{L^{2}(\mathbb{T})}. \label{eq2.18}%
\end{equation}
(Condition (\ref{eq2.12}) is not needed for this, but only for the algorithmic
operations of the Appendices in Refs.\ \citenum{BrJo99b} and \citenum{Jor00}.)
We also showed that the solutions $(m_{j}%
)_{j=0}^{N-1}$ to (\ref{eq2.15}) are in $1$--$1$ correspondence with the group
of all polynomial functions
\begin{equation}
A\colon\mathbb{T}\longrightarrow\mathrm{U}_{N}(\mathbb{C}) \label{eq2.19}%
\end{equation}
where $\mathrm{U}_{N}(\mathbb{C})$ denotes the unitary $N\times N$ matrices.
The correspondence is $m\leftrightarrow A$ with
\begin{equation}
m_{j}(z)=\sum_{k=0}^{N-1}A_{j,k}(z^{N})z^{k}, \label{eq2.20}%
\end{equation}
and%
\begin{equation}
A_{j,k}(z)=\frac{1}{N}\sum_{w^{N}=z}w^{-k}m_{j}(w). \label{eq2.21}%
\end{equation}
This correspondence plays a central role in the proofs
in Refs.\ \citenum{BrJo99b} and \citenum{Jor00}. We also show
in Ref.\ \citenum{BrJo00} that if $m_{0}$ is given, 
and if it satisfies (\ref{eq2.11}),
then it is possible to construct $m_{1},\dots,m_{N-1}$ such that the extended
system $m_{0},m_{1},\dots,m_{N-1}$ will satisfy (\ref{eq2.15}). As a
consequence, $A$ in (\ref{eq2.21}) will be a $\mathrm{U}_{N}(\mathbb{C}%
)$-loop, i.e., $A\colon\mathbb{T}\rightarrow\mathrm{U}_{N}(\mathbb{C})$, and
moreover, the original $m_{0}$ is then recovered from (\ref{eq2.20}) for $j=0$.

By virtue of (\ref{eq2.17})--(\ref{eq2.18}), $L^{2}(\mathbb{T})$, or
equivalently $\ell^{2}(\mathbb{Z})$, splits up as an orthogonal sum%
\begin{equation}
S_{j}(\ell^{2}(\mathbb{Z}))\text{,\qquad}j=0,1,\dots,N-1.\label{eq2.22}%
\end{equation}
We saw that the wavelet transform $W$ of (\ref{eq2.9})--(\ref{eq2.10}) maps
$\ell^{2}(\mathbb{Z})$ onto $\mathcal{V}_{0}$, and from (\ref{eq2.13}) we
conclude that $W$ maps $S_{0}(\ell^{2}(\mathbb{Z}))$ onto $U_{N}%
(\mathcal{V}_{0})$ ($=:\mathcal{V}_{1}$). Hence, in the $N$-scale wavelet
case, $W$ transforms the spaces $S_{j}(\ell^{2}(\mathbb{Z}))$ ($\subset
\ell^{2}(\mathbb{Z})$) onto orthogonal subspaces $\mathcal{W}_{1}^{(j)}$,
$j=1,\dots,N-1$ in $L^{2}(\mathbb{R})$, and
\begin{equation}
\mathcal{W}_{1}=\mathcal{V}_{0}\ominus\mathcal{V}_{1}=\sum_{j=1}%
^{N-1}\!\vphantom{\sum}^{\raisebox{3pt}[0pt][0pt]{$\scriptstyle\oplus$}%
}\,\mathcal{W}_{1}^{\left(  j\right)  },\label{eq2.25}%
\end{equation}
where
\begin{equation}
\mathcal{W}_{1}^{\left(  j\right)  }=S_{j}\ell^{2},\qquad j=1,\dots
,N-1.\label{eq2.25bis}%
\end{equation}
Each of the spaces $\mathcal{V}_{1}$ and $\mathcal{W}_{1}^{\left(  j\right)
}$ is split further into orthogonal subspaces corresponding to iteration of
the operators $S_{0},S_{1},\dots,S_{N-1}$ of (\ref{eq2.17})--(\ref{eq2.18})$.$
It is the system $\{S_{j}\}_{j=0}^{N-1}$ which is called a wavelet
representation, and it follows that the wavelet decomposition may be recovered
from the representation. Moreover, the variety of all wavelet representations
is in $1$--$1$ correspondence with the group of polynomial functions $A$ as
given in (\ref{eq2.19}). The correspondence is fixed by (\ref{eq2.20}%
)--(\ref{eq2.21}). Operators $\{S_{j}\}$ satisfying (\ref{eq2.17}%
)--(\ref{eq2.18}) are said to constitute a representation of the $C^{\ast}%
$-algebra $\mathcal{O}_{N}$, the Cuntz algebra\cite{Cun77}, and it is the
irreducibility of the representations from (\ref{eq2.16}) which will concern
us. If a representation (\ref{eq2.16}) is reducible (Definition \ref{DefRep.3}%
), then there is a subspace
\begin{equation}
0\subsetneqq\mathcal{H}_{0}\subsetneqq L^{2}(\mathbb{T})\label{eq2.24}%
\end{equation}
which is invariant under all the operators $S_{j}$ and $S_{j}^{\ast}$, and so
the data going into the wavelet filter system $\{m_{j}\}$ are then not minimal.

\section{The main theorem}

\label{The}

The main result will be stated in the present section, but without proof.
Instead the reader is referred to Ref.\ \citenum{Jor00} 
for the full proof, and for a
detailed discussion of its implications. We noted above that the
representation (\ref{eq2.16}) given from a QMF system $m_{j}=m_{j}^{\left(
A\right)  }$ via (\ref{eq2.20})--(\ref{eq2.21}) is irreducible if and only if
the subbands are optimal, in that they do not admit further reduction into a
refined system of closed subspaces of $L^{2}\left(  \mathbb{R}\right)  $.

\begin{theorem}
\label{ThmThe}The representation
\[
S_{j}^{\left(  A\right)  }f\left(  z\right)  =m_{j}^{\left(  A\right)
}\left(  z\right)  f\left(  z^{N}\right)  ,\qquad f\in L^{2}\left(
\mathbb{T}\right)  ,\;z\in\mathbb{T}%
\]
is an irreducible representation of $\mathcal{O}_{N}$ on $L^{2}\left(
\mathbb{T}\right)  $ if and only if $A\colon\mathbb{T}\rightarrow
\mathrm{U}_{N}\left(  \mathbb{C}\right)  $ does not admit a matrix corner of
the form
\begin{equation}
V%
\begin{pmatrix}
z^{n_{0}} &  &  & \\
& z^{n_{1}} &  & \llap{\smash{\huge$0$}}\\
&  & \ddots & \\
\rlap{\smash{\huge$0$}} &  &  & z^{n_{M-1}}%
\end{pmatrix}
,\label{eqThmMin.7AzVm}%
\end{equation}
for some $V\in\mathrm{U}_{M}\left(  \mathbb{C}\right)  $, and where
$n_{0},n_{1},\dots,n_{M-1}\in\left\{  0,1,2,\dots\right\}  $.
\end{theorem}

\section*{ACKNOWLEDGMENTS}
We gratefully acknowledge the
expert \TeX-typesetting and
manuscript preparation by Brian Treadway,
as well as partial support by the National Science Foundation.

\bibliographystyle{spiebib}
\bibliography{jorgen}
\end{document}